\documentclass[12pt]{article}
\usepackage{amssymb,amsmath,amsthm,bm}
\usepackage{color}
\newcommand{\lcolor}[1]{\textcolor{black}{#1}}  

\begin{document}
\newcommand{\bea}{\begin{eqnarray}}
\newcommand{\ena}{\end{eqnarray}}
\newcommand{\beas}{\begin{eqnarray*}}
\newcommand{\enas}{\end{eqnarray*}}
\newcommand{\beq}{\begin{equation}}
\newcommand{\enq}{\end{equation}}
\def\qed{\hfill \mbox{\rule{0.5em}{0.5em}}}
\newcommand{\ignore}[1]{}
\newcommand{\transpose}{{\mbox{\scriptsize\sf T}}}
\newcommand{\N}{{\bf n}}
\newcommand{\spa}{S}
\newcommand{\bbs}{\mathbb{S}}
\newcommand{\PP}{P}
\newcommand{\vecS}{{\bf U}}
\newcommand{\vecR}{{\bf V}}
\newcommand{\qmq}[1]{\quad \mbox{#1} \quad}
\newcommand{\aBKR}{(a)}
\newcommand{\bBKR}{(b)}
\newcommand{\cBKR}{(c)}
\newcommand{\dBKR}{(d)}
\newcommand{\eBKR}{(e)}
\newcommand{\apBKR}{(a$\,'$) }
\newcommand{\bpBKR}{(b$\,'$) }
\newcommand{\cpBKR}{(c$\,'$) }
\newcommand{\dpBKR}{(d$\,'$) }
\newcommand{\epBKR}{(e$\,'$) }
\newtheorem{theorem}{Theorem}[section]
\newtheorem{corollary}{Corollary}[section]
\newtheorem{conjecture}{Conjecture}[section]
\newtheorem{proposition}{Proposition}[section]
\newtheorem{lemma}{Lemma}[section]
\newtheorem{definition}{Definition}[section]
\newtheorem{example}{Example}[section]
\newtheorem{framework}{Framework}
\newtheorem{remark}{Remark}[section]
\newtheorem{case}{Case}[section]
\newtheorem{condition}{Condition}[section]
\newcounter{saveenum}

\def\undertilde#1{\mathord{\vtop{\ialign{##\crcr
				$\hfil\displaystyle{#1}\hfil$\crcr\noalign{\kern1.5pt\nointerlineskip}
				$\hfil\tilde{}\hfil$\crcr\noalign{\kern1.5pt}}}}}

\title{Functional van den Berg-Kesten-Reimer Inequalities and their Duals, with Applications}

\author{Larry Goldstein and Yosef Rinott\thanks{Research supported by Israel Science Foundation Grant No.
473/04}\\University of Southern California, Hebrew University of
Jerusalem} 
\maketitle

\date{}

\begin{center}  Abbreviated Title: Functional BKR Inequalities
\end{center}

\begin{abstract}
The BKR inequality conjectured by van den Berg and Kesten in [\ref{vbk}], and proved by Reimer in [\ref{reim}], states that for
$A$ and $B$ events on $S$, a finite product of finite sets $S_i,i=1,\ldots,n$, and $P$ any product measure on $S$, 
$$ P(A \Box B) \le P(A)P(B),$$ where the set $A \Box B$ consists of the elementary events which lie in both $A$ and $B$ for `disjoint reasons.' 
Precisely, with ${\bf n}:=\{1,\ldots,n\}$ and $K \subset {\bf n}$, for ${\bf x} \in S$ letting $[{\bf x}]_K=\{{\bf y} \in S: y_i = x_i, i \in K\}$, 
the set $A \Box B$ consists of all ${\bf x} \in S$ for which there exist disjoint subsets $K$ and $L$ of ${\bf n}$ for which $[{\bf x}]_K \subset A$ and $[{\bf x}]_L \subset B$.

The BKR inequality is extended to the following functional version on a general finite product measure space $(S,\mathbb{S})$  with product probability measure $P$, 
$$E\left\{
\max_{\stackrel{K \cap L = \emptyset}{K \subset {\bf n}, L \subset {\bf n}}}
		\underline{f}_K({\bf X})\underline{g}_L({\bf
	X})\right\} \leq E\left\{f({\bf X})\right\}\,E\left\{g({\bf
	X})\right\},$$
where $f$ and $g$ are non-negative measurable functions,
$\underline{f}_K({\bf x}) = {\rm ess} \inf_{{\bf y}
\in [{\bf x}]_K}f({\bf y})$ and
$\underline{g}_L({\bf x}) = {\rm ess} \inf_{{\bf y} \in [{\bf
x}]_L}g({\bf y}).$ The original BKR inequality is recovered  by taking $f({\bf x})={\bf 1}_A({\bf x})$ and $g({\bf x})={\bf 1}_B({\bf x})$,
and applying the fact that in general ${\bf 1}_{A \Box B} \le \max_{K \cap L = \emptyset}
\underline{f}_K({\bf x}) \underline{g}_L({\bf x})$.

Related formulations, and functional versions of the dual
inequality on events by Kahn, Saks, and Smyth [\ref{kss}], are
also considered. Applications include order statistics, assignment
problems, and paths in random graphs.
\end{abstract}


\section{Introduction}
\lcolor{This paper is a minor revision of [\ref{GoRi}], correcting an error in equation \eqref{was.5} that was pointed out to us by Richard Arratia during the preparation of [\ref{AGH}], a draft of which he shared with us. Inequality \eqref{was.5} was earlier incorrectly stated as an equality. While the correction is important on its own, the error is inconsequential for the purposes of our original work. As \eqref{was.5} was only applied to show \eqref{boximpliesbkrforindicators}, for which we now provide a much simpler argument not involving \eqref{was.5}, no results depended on its validity. In particular, the statements of all theorems here are the same as in [\ref{GoRi}].}

For ${\bf x} =(x_1,\ldots,x_n) \in \spa$, where $\spa=\prod_{i=1}^n
\spa_i$ any product space, and $K =\{k_1, \ldots,k_m\} \subseteq
\N:=\{1,\ldots,n\}$ with $k_1<\cdots<k_m$, define
$${\bf x}_K=(x_{k_1},\ldots,x_{k_m}) \quad {\rm and} \quad [{\bf
x}]_K = \{ {\bf y} \in \spa: y_K=x_K \},
$$ the restriction of ${\bf x}$ to the indicated coordinates,
and the collection of all elements in $\spa$ which agree with
${\bf x}$ in those coordinates, respectively. For $A,B \subseteq
\spa$  we say that ${\bf x} \in A,{\bf y} \in B$ {\em disjointly}
if there exist \bea \label{def-disjoint} \mbox{$K,L \subseteq \N,
K \cap L = \emptyset$ such that} \quad [{\bf x}]_K \subseteq A
\quad \mbox{and} \quad [{\bf y}]_L \subseteq B, \ena and denote
\bea \label{gen-box}
A \Box B = \{{\bf x}: {\bf x} \in A, {\bf x}
\in B \,\mbox{disjointly} \}.
\ena
The operation $A \Box B$ corresponds to elementary events
which are in both $A$ and $B$ for disjoint `reasons' in the sense
that inclusion in $A$ and $B$ is determined on disjoint sets of
coordinates.

Theorem \ref{BKR} was conjectured in van den Berg and Kesten
[\ref{vbk}]. It was proved in [\ref{vbk}] for $A$ and $B$
increasing sets and $\spa=\{0,1\}^n$, and it was also demonstrated
there that Theorem \ref{BKR} follows from its special case
$\spa=\{0,1\}^n$. Using the latter fact, the conjecture was
established in general by Reimer [\ref{reim}].

\begin{theorem}
\label{BKR} For $P=\prod_{i=1}^n P_i$ any product measure on
$\spa=\prod_{i=1}^n \spa_i$, $\spa_i$ finite,
\begin{equation}
\label{original-box} P(A \Box B) \leq P(A)P(B).
\end{equation}
\end{theorem}

Many useful formulations can be found in van den Berg and Fiebig
[\ref{vbf}], in addition to the following motivating example which
appeared earlier in [\ref{vbk}]. Independently assign a random
direction to each edge $e=\{v_i,v_j\}$ of a finite graph, with
$p_e(v_i,v_j)=1-p_e(v_j,v_i)$ the probability of the edge $e$
being directed from vertex $v_i$ to $v_j$. With $V_1,V_2,W_1,W_2$
sets of vertices, Theorem \ref{BKR} yields that the product of the
probabilities that there exist directed paths from $V_1$ to $V_2$
(event $A$) and from $W_1$ to $W_2$ (event $B$) is an upper bound
to the probability that there exist two disjoint directed paths,
one from $V_1$ to $V_2$ and another from $W_1$ to $W_2$ (event $A
\Box B$).

The main thrust of this paper is to show how Theorem \ref{BKR}
implies inequalities in terms of functions, of which
(\ref{original-box}) is the special case of indicators, and
similarly for the dual inequality of [\ref{kss}]. These functional
inequalities, and their duals, are stated in Theorems \ref{thm1} and
\ref{thm1-dual}, and their proofs can be found in Section
\ref{proof}. Applications to order statistics, allocation problems,
and random graphs are given in Section \ref{appl}. Specializing to
monotone functions, we derive related inequalities and stochastic
orderings in Section \ref{PQD-section}; these latter results are
connected to those of Alexander [\ref{alex}].

For each $i=1,\ldots,n$, let $(\spa_i,\bbs_i)$ be measurable
spaces, and set $\spa=\prod_{i=1}^n \spa_i$ and
$\bbs=\bigotimes_{i=1}^n \bbs_i$, the product sigma algebra.
Henceforth, all given real valued functions on $\spa$, such as
$f_\alpha,g_\beta, \alpha \in {\cal A}, \beta \in {\cal B}$ are
assumed to be $(\bbs,\mathbb{B})$ measurable where $\mathbb{B}$
denotes the Borel sigma algebra of ${\bf R}$, and functions on $S$
with values in $2^{\N}$, such as $K({\bf x})$ in inequality \dBKR
\,\, of Theorem \ref{thm1} below, are assumed to be
$(\bbs,2^{2^{\N}})$ measurable. Measurability issues arise in
definitions (\ref{Kalphax}), (\ref{fgismin}), and
(\ref{maxKL-measurable}), and are settled in Section
\ref{appendix}. We also show in Section \ref{appendix} that
Theorem \ref{thm1} applies to the completion of the measure space
$(\spa,\bbs)$ with respect to the measure $P$ appearing in the
theorem; similarly for Theorem \ref{thm1-dual}.

For $K \subseteq \N$ we say that a function $f$ defined on $\spa$
depends on $K$ if ${\bf x}_K={\bf y}_K$ implies $f({\bf x})=f({\bf
y})$. The inequalities in Theorems \ref{thm1} and \ref{thm1-dual}
require one of two frameworks, the first of which is the
following.
\begin{framework}
\label{I} $\{f_\alpha({\bf x})\}_{ \alpha \in {\cal A}}$\,\, and
\,\,$\{g_\beta({\bf y})\}_{\beta \in {\cal B}}$ are given
collections of non-negative functions on $\spa$,  such that
$f_\alpha, g_\beta$ depend respectively on subsets of $\N$
$K_\alpha, L_\beta$ in ${\cal K}=\{K_\alpha\}_{\alpha \in {\cal
A}}$ and ${\cal L}=\{L_\beta\}_{\beta \in {\cal B}}$, where ${\cal
A}$ and ${\cal B}$ are finite or countable.
\end{framework}

The elements of ${\cal K}$ and ${\cal L}$ are not assumed to be
distinct; we may have, say, $K_\alpha=K_{\gamma}$ for some $\alpha
\ne \gamma$ and $f_\alpha \ne f_{\gamma}$.  Note also that if a
function depends on $K$, it depends on any subset of $\N$
containing $K$. For notational brevity we may write $\alpha$ for
$K_\alpha$; for example, we may use $\alpha \cap \beta$ as an
abbreviation for $K_\alpha \cap L_\beta$, and also ${\bf
x}_\alpha$ for ${\bf x}_{K_\alpha}$.

The second framework is
\begin{framework}
\label{II} $f$ and $g$ are two given non-negative functions, and
$\cal K$ and $\cal L$ are any subsets of $2^{\bf n}$. With
$P$ a probability measure on $(S,\mathbb{S})$ define for $K \in
{\cal K}, L \in {\cal L}$,
\begin{equation}
\label{fgismin} \underline{f}_K({\bf x}) = {\rm ess} \inf_{{\bf y}
\in [{\bf x}]_K}f({\bf y}), \quad \mbox{and} \quad
\underline{g}_L({\bf x}) = {\rm ess} \inf_{{\bf y} \in [{\bf
x}]_L}g({\bf y}),
\end{equation}
where the essential infimums for $\underline{f}_K({\bf x})$ and $\underline{g}_L({\bf x})$ are being taken with respect to the product probability measure on the coordinates in $K^c$ and $L^c$ respectively.
\end{framework}

Our functional extension of the BKR inequality
(\ref{original-box}) is
\begin{theorem}
\label{thm1} Let ${\bf X} = (X_1,\ldots,X_n) \in \spa$ be a random
vector and $P$ a probability measure on $(S,\mathbb{S})$ under
which $X_1,\ldots,X_n$ are independent.
\begin{enumerate}
\item Under framework \ref{I},
$$
E\left\{\sup_{ \alpha \cap \beta = \emptyset} f_\alpha({\bf
X})g_\beta({\bf X})\right\} \leq E\left\{\sup_\alpha f_\alpha({\bf
X})\right\} \,E\left\{\sup_\beta g_\beta({\bf X})\right\}.
\eqno{(a)}
$$

\item Under framework \ref{II},
$$
E\left\{\max_{\stackrel{K \cap L = \emptyset}{K \in {\cal K}, L
\in {\cal L}}}\underline{f}_K({\bf X})\underline{g}_L({\bf
X})\right\} \leq E\left\{f({\bf X})\right\}\,E\left\{g({\bf
X})\right\}. \eqno{(b)}
$$
\end{enumerate}
\end{theorem}
The special case of \bBKR\,where ${\cal K}={\cal L}=2^{\bf n}$, the
collection of all subsets of ${\bf n}$, clearly implies the
inequality in general.

In [\ref{reim}], inequality \eqref{original-box} for the $\Box$ operation was proven only for discrete finite product spaces, that is, a finite product of finite sets; Theorem \ref{thm1} applies to functions of a vector ${\bf X}$ having independent coordinates taking values in \textit{any} measure
space. \lcolor{For $f({\bf x})={\bf 1}_A({\bf x})$ and $g({\bf x})={\bf 1}_B({\bf x})$ for $A,B \in {\mathbb S}$, we have
\bea \label{was.5}
{\bf 1}_{A \Box B} \le \max_{K \cap L = \emptyset}
\underline{f}_K({\bf x}) \underline{g}_L({\bf x}).
\ena
To see \eqref{was.5}, note that replacing essential infimum by infimum in \eqref{fgismin}, the inequality becomes equality. Hence \eqref{was.5} holds as stated because  the essential infimum is at least as large as the infimum. In other words, elements of $A \Box B$ demand disjoint `reasons' for $A$ and $B$ that hold for all outcomes in the probability space, while the right hand side of \eqref{was.5} only requires that the `reasons' be almost sure.}

\lcolor{
In [\ref{GoRi}], the specialization of \bBKR\, to indicator functions and ${\cal K}={\cal L}=2^{\bf n}$, and `equality' in \eqref{was.5}, was interpreted to mean that inequality \eqref{original-box} holds for general product spaces. However, as \bBKR\, is an inequality, this interpretation now yields that for all $A,B \in {\mathbb S}$
the set in $\mathbb{S}$ whose indicator appears in the right hand side of \eqref{was.5}
contains $A \Box B$ and has probability bounded above by the product of the probabilities of $A$ and $B$.
This latter interpretation appears in
Corollary 4 of [\ref{AGH}] when $\mathbb{S}$ is taken to be Euclidean space.}

The following is a straightforward generalization of Theorem
\ref{thm1}, stated here only for inequality \aBKR. Note that in the inequality below, as $m$ increases the pairwise constraints
$\alpha_i \cap \alpha_j = \emptyset$  become more restrictive, and
the inequality less sharp.

\begin{theorem}
Let ${\bf X} \in \spa$ be a random vector with independent
coordinates. Then for given finite or countable collections of
non-negative functions $\{f_{i,\alpha}\}_{\alpha \in {\cal A}_i}$
depending on $\{K_{i,\alpha}\}_{\alpha \in {\cal A}_i}$\,,
\,\,$i=1, \ldots, m,\,$
\begin{equation*} 
E\left\{\sup_{ \stackrel{(\alpha_1,\ldots,\alpha_m) \in {\cal A}_1
\times \cdots \times {\cal A}_m}{\alpha_k \cap \alpha_l =
\emptyset, \,\, k \not = l}}\,\, \prod_{i=1}^m f_{i,\alpha_i}({\bf
X})\right\} \leq \prod_{i=1}^m E\left\{\sup_{\alpha \in {\cal
A}_i} f_{i, \alpha}({\bf X})\right\}.
\end{equation*}
\end{theorem}

Next we describe an inequality of Kahn, Saks, and Smyth
[\ref{kss}], which may be considered dual to the BKR inequality
(\ref{original-box}), and then provide a function version. We use
a notation compatible with (\ref{original-box}). With `disjointly'
defined in (\ref{def-disjoint}), denote \beas
A \Diamond B = \{({\bf x},{\bf y}): {\bf x}
\in A, {\bf y} \in B \,\mbox{disjointly}\}. \enas
Note that
$$
{\bf 1}_{A \Box B}({\bf x})={\bf 1}_{A \Diamond B}({\bf x},{\bf
	x}).
$$

The following,
which we call the KSS inequality, is dual to Theorem \ref{BKR} and
is given in [\ref{kss}].
\begin{theorem}
\label{kss-dual-reimer} If $P$ denotes the uniform measure over
$\{0,1\}^n$, then for  any $(A,B) \subseteq \{0,1\}^n \times
\{0,1\}^n$, \bea \label{dual-reimer} (P\times P)(A \Diamond B) \le
P(A \cap B). \ena
\end{theorem}

Our functional extension of the KSS inequality is as follows.
\begin{theorem}
\label{thm1-dual} Let ${\bf X} = (X_1,\ldots,X_n) \in \spa$ be a
random vector, $P$ any probability measure on $(S,\mathbb{S})$
such that $X_1,\ldots,X_n$ are independent, and ${\bf Y}$ an
independent copy of $\bf X$.
\begin{enumerate}
\item
Under framework \ref{I},
$$
E \left\{\sup_{\alpha \cap \beta = \emptyset}f_\alpha({\bf
X})g_\beta({\bf Y})\right\} \le E
\left\{\sup_{\alpha,\beta}f_\alpha({\bf X})g_\beta({\bf
X})\right\}. \eqno{(a')}
$$

\item Under framework \ref{II},
$$
E\left\{\max_{\stackrel{K \cap L = \emptyset}{K \in {\cal K},L \in
{\cal L}}}\underline{f}_K({\bf X})\underline{g}_L({\bf Y})\right\}
\leq E\left\{f({\bf X}) g({\bf X})\right\}. \eqno{(b')}
$$
\end{enumerate}
\end{theorem}

The $\Diamond$ operation in Theorem \ref{kss-dual-reimer} on
$\{0,1\}^n \times \{0,1\}^n$ was defined implicitly in [\ref{kss}],
and the inequality was extended there to product measure on discrete
finite product spaces. With $f({\bf x})$ and $g({\bf x})$ the
indicator functions of sets $A$ and $B$ respectively, \lcolor{we have the inequality} \bea \label{DnBox} {\bf 1}_{A \Diamond B}({\bf
x},{\bf y}) \lcolor{\le} \max_{K \cap L = \emptyset} \underline{f}_K({\bf x})
 \underline{g}_L({\bf y}),
\ena where $\underline{f}_K,\underline{g}_L$ are given in
(\ref{fgismin}). Therefore, inequality \bpBKR of Theorem
\ref{thm1-dual} specialized to the case where ${\cal K}={\cal
L}=2^{\bf n}$ and $f$ and $g$ are indicators says that the
original KSS inequality \lcolor{\eqref{dual-reimer}} for events in discrete finite product
spaces extends to vectors having independent coordinates taking values in any measure space \lcolor{in the sense that $A \Diamond B$ is a subset of the set whose indicator is the function appearing on the right hand side of \eqref{DnBox}, and has probability bounded by $P(A \cap B)$, by (b').}

We next discuss further formulations of Theorems \ref{thm1} and
\ref{thm1-dual} which are of independent interest, and will be
used in the proof. Under Framework \ref{I}, for any subsets $K$ and $L$
of $\N$, define
\begin{equation}
\label{one-to-2} \tilde{f}_K({\bf x})=\sup_{\alpha:K_\alpha
\subseteq K}f_\alpha({\bf x}) \quad \mbox{and} \quad
\tilde{g}_L({\bf x})=\sup_{\beta:L_\beta \subseteq L}g_\beta({\bf
x}).
\end{equation}
For any given functions $K({\bf x})$ and $L({\bf x})$ defined on
$\spa$ and taking values in $2^{\bf n}$, under Framework \ref{I},
extend (\ref{one-to-2}) to \bea \label{Kalphax} \tilde{f}_{K({\bf
x})}({\bf x}) = \sup_{\alpha: K_\alpha \subseteq K({\bf x})}
f_\alpha({\bf x}) \quad \mbox{and} \quad \tilde{g}_{L({\bf x})}({\bf
x}) = \sup_{\beta: L_\beta \subseteq L({\bf x})} g_\beta({\bf x}),
\ena and under Framework \ref{II}, extend (\ref{fgismin}) to
\begin{equation}
\label{fgisminx} \underline{f}_{K({\bf x})}({\bf x}) = {\rm ess}
\inf_{{\bf y} \in [{\bf x}]_{K({\bf x})}}f({\bf y}), \quad
\mbox{and} \quad \underline{g}_{L({\bf x})}({\bf x}) = {\rm ess}
\inf_{{\bf y} \in [{\bf x}]_{L({\bf x})}}g({\bf y}).
\end{equation}
Proposition \ref{reformulations} shows that parts \cBKR\, and \dBKR\, of
Proposition \ref{prop1} are reformulations of \aBKR\, of Theorem
\ref{thm1}, and likewise \eBKR\, a reformulation of \bBKR.

\begin{proposition}
\label{prop1} Let the hypotheses of Theorem \ref{thm1} hold. In
Framework \ref{I} \cBKR\, and \dBKR \,\, below obtain.
$$
E\left\{ \max_{\stackrel{K \cap L = \emptyset}{K \in {\cal K}, L
\in {\cal L}}} \tilde{f} _{K}({\bf X}) \tilde{g}_{L}({\bf X})
\right\} \leq E\left\{\sup_\alpha f_\alpha({\bf X}) \right\}
\,E\left\{\sup_\beta g_\beta({\bf X})\right\}. \eqno{(c)}
$$
$$
E\left\{\tilde{f} _{K({\bf X})}({\bf X}) \tilde{g}_{L({\bf
X})}({\bf X})\right\} \leq E\left\{\sup_\alpha f_\alpha({\bf
X})\right\} \,E\left\{\sup_\beta g_\beta({\bf X})\right\},
\eqno{(d)}
$$
holding for any given  $K({\bf x}) \in {\cal K}$ and $L({\bf x})
\in {\cal L}$ such that $K({\bf x}) \cap L({\bf x})=\emptyset$.
\smallskip
\smallskip

In Framework \ref{II} we have
$$
E\left\{ \underline{f}_{K({\bf X})}({\bf X})
\underline{g}_{{L}({\bf X})}({\bf X})\right\} \leq E\left\{f({\bf
X})\}\,E\{g({\bf X})\right\}, \eqno{(e)}
$$
holding for any given  $K({\bf x}) \in {\cal K}$ and $L({\bf x})
\in {\cal L}$ such that $K({\bf x}) \cap L({\bf x})=\emptyset$.
\end{proposition}

As in Theorem \ref{thm1} the special cases of \cBKR\, and \dBKR\,
where ${\cal K}$ and ${\cal L}$ both equal $2^{\bf n}$ implies the
inequality in general. Similarly, the special case of inequality
\eBKR\, with ${\cal K}= {\cal L} = 2^\N$ and $L({\bf x})=K^c({\bf
x})$, where $K^c$ denotes the complement of $K$, yields the
inequality in general, that is, \eBKR\, is equivalent to the
statement that for any given $K({\bf x})$,
$$E\left\{ \underline{f}_{K({\bf X})}({\bf X})
\underline{g}_{{K^c}({\bf X})}({\bf X})\right\} \leq E\left\{f({\bf
X})\}\,E\{g({\bf X})\right\}.$$

Parallel to the claims of Proposition \ref{prop1}, parts \cpBKR and
\dpBKR below are reformulations of \apBKR of Theorem
\ref{thm1-dual}, and \epBKR a reformulation of \bpBKR. For given $f$
and a function $K({\bf x},{\bf y})$ taking values in $2^{\bf n}$,
define $\tilde{f} _{K({\bf X},{\bf Y})}$ and $\underline{f}_{K({\bf
X},{\bf Y})}$ by replacing $K({\bf x})$ by $K({\bf x},{\bf y})$ in
(\ref{Kalphax}) and (\ref{fgisminx}) respectively.

\begin{proposition}
\label{prop1-dual} Let the hypotheses of Theorem \ref{thm1-dual}
hold. In Framework \ref{I} \cpBKR and \dpBKR below obtain.
$$E\left\{\max_{K \in {\cal K},L \in {\cal L} }\tilde{f} _{K}({\bf
X}) \tilde{g}_{L}({\bf Y})\right\} \leq E\left\{\sup_{\alpha \cap
\beta=\emptyset} f_\alpha({\bf X}) g_\beta({\bf X})\right\}.
\eqno{(c\,')}
$$

$$
E\left\{\tilde{f} _{K({\bf X},{\bf Y})}({\bf X}) \tilde{g}_{L({\bf
X},{\bf Y})}({\bf Y})\right\} \leq E\left\{\sup_{\alpha \cap \beta
= \emptyset} f_\alpha({\bf X}) g_\beta({\bf X})\right\},
\eqno{(d\,')}
$$
holding for any given $K({\bf x} ,{\bf y})\in {\cal K}$ and
$L({\bf x},{\bf y})\in {\cal L}$ replacing $K{(\bf x})$ and
$L({\bf x})$ in (\ref{Kalphax}), respectively, and satisfying
$K({\bf x},{\bf y}) \cap L({\bf x},{\bf y})= \emptyset$.

\smallskip

In Framework \ref{II} we have

$$E\left\{ \underline{f}_{K({\bf X},{\bf Y})}({\bf X})
\underline{g}_{{L}({\bf X},{\bf Y})}({\bf Y})\right\} \leq
E\left\{f({\bf X})g({\bf X})\right\}, \eqno{(e\,')}
$$
for given $K({\bf x},{\bf y})  \in {\cal K}$ and $L({\bf x},{\bf
y}) \in {\cal L}$ replacing $K{(\bf x})$ and $L({\bf x})$ in
(\ref{fgisminx}), respectively.
\end{proposition}

\section{Applications}
\label{appl}

\begin{example} \textbf{Order Statistics Type Inequalities}
\label{order-statistics} Let ${\bf X}=(X_1,\ldots,X_n)$ be a
vector of independent non-negative random variables with
associated order statistics $X_{[n]} \le \cdots \le X_{[1]}$. Let
${\cal A}={\cal B}$ be the collection of all the singletons
$\alpha \in \N$ and $f_\alpha({\bf x})=g_\alpha({\bf
x})=x_\alpha$. Then
$$
\max_\alpha f_\alpha({\bf X})=X_{[1]}, \quad \max_{\alpha \cap
\beta = \emptyset}f_\alpha({\bf X}) g_\beta({\bf X}) =
X_{[1]}X_{[2]},
$$
and inequality \aBKR\, of Theorem \ref{thm1} provides the middle
inequality in the string
\beas
EX_{[1]}EX_{[2]} \le EX_{[1]} X_{[2]} \le (EX_{[1]})^2 \le
EX_{[1]}^2.
\enas
The leftmost inequality is true since order
statistics are always positively correlated (moreover they are
\textit{associated} as defined by Esary et al [\ref{epw}], and
even MTP$_2$, see Karlin and Rinott [\ref{kr}]); the rightmost
inequality follows from Jensen.

Theorem \ref{thm1} allows a large variety of extensions of this
basic order statistics inequality. For example, taking ${\cal A}$
and ${\cal B}$ to be all $k$ and $l$ subsets of $\N$ respectively,
then with \bea \label{fgprod} f_\alpha({\bf x})=\prod_{j \in
\alpha}x_j \quad \mbox{and} \quad g_\beta({\bf x})=\prod_{j \in
\beta}x_j \ena we derive
$$
E\left( \prod_{j=1}^{k+l}X_{[j]} \right) \le E\left( \prod_{j=1}^k
X_{[j]}\right) E\left( \prod_{j=1}^l X_{[j]}\right).
$$

Dropping the non-negativity assumption on $X_1,\ldots,X_n$, we
have for all $t>0$,
$$Ee^{t(X_{[1]}+X_{[2]})} \le
[Ee^{tX_{[1]}}]^2=Ee^{tX_{[1]}}Ee^{tY_{[1]}}=Ee^{t(X_{[1]}+Y_{[1]})},$$
with $Y_i$'s being independent copies of the $X_i$'s. Likewise,
for all $t>0$,
$$Ee^{-t(X_{[n]}+X_{[n-1]})} \le
[Ee^{-tX_{[n]}}]^2=Ee^{-tX_{[n]}}Ee^{-tY_{[n]}}=Ee^{-t(X_{[n]}+Y_{[n]})}.$$
Moment generating function and Laplace orders are discussed in
Shaked and Shanthikumar [\ref{ss}].

Returning to non-negative variables, a  variation  of
(\ref{fgprod}) follows by replacing products with sums, that is,
\bea \label{fgsum} f_\alpha({\bf x})=\sum_{j \in \alpha}x_j \quad
\mbox{and} \quad g_\beta({\bf x})=\sum_{j \in \beta}x_j, \ena
which for, $k=l=2$ say, yields
$$
E\max_{\{i,j,k,l\}=\{1,2,3,4\}}(X_{[i]}+X_{[j]})(X_{[k]}+X_{[l]})
\le [E(X_{[1]}+X_{[2]})]^2.
$$
Though the maximizing indices on the left hand side will be
$\{1,2,3,4\}$ as indicated, the choice is not fixed and depends on
the $X$'s; note, for example, that
$(X_{[1]}+X_{[2]})(X_{[3]}+X_{[4]})$ is never maximal apart from
degenerate cases.

Definition  (\ref{fgprod}) and (\ref{fgsum}) are special cases
where $f$ and $g$ are increasing non-negative functions of $k$ and
$l$ variables and \bea \label{fixedfgs} f_\alpha({\bf x}) = f({\bf
x}_\alpha) \quad \mbox{and} \quad g_\beta({\bf x}) = g({\bf
x}_\beta); \ena when $f$ and $g$ are symmetric,
\begin{multline*}
E\max_{\{i_1,\ldots,i_k,j_1,\ldots,j_{\,l\}}=\{1,\ldots,k+l\}}f(X_{[i_1]},\ldots,X_{[i_k]})g
(X_{[j_1]},\ldots,X_{[j_{\hspace{.01in}l}]})\\ \le
Ef(X_{[1]},\ldots,X_{[k]}) Eg(X_{[1]},\ldots,X_{[l]}).
\end{multline*}

We now give an example which demonstrates that these order
statistics type inequalities can be considered in higher
dimensions. Let ${\bf X}_1,\ldots,{\bf X}_n$ be independent
vectors in ${\bf R}^m$, and for $\alpha, \beta \subseteq \N$ with
$|\alpha|=|\beta|=3$ let $f_\alpha$ and $g_\beta$ be given as in
(\ref{fixedfgs}), where $f({\bf x}_1,{\bf x}_2,{\bf x}_3)=g({\bf
x}_1,{\bf x}_2,{\bf x}_3)$ is, say, the area of the triangle
formed by the given three vectors.  Theorem \ref{thm1} gives that
the expected greatest product of the areas of two triangles with
distinct vertices is bounded above by the square of the
expectation of the largest triangular area.

To explore the dual inequality in these settings, let ${\bf X}$ be
a vector of independent variables with support contained in
$[0,1]$, and $\bf Y$ an independent copy.  With  ${\cal A}={\cal
B}$  the collections of all singletons $\alpha$ in $\N$, and
$f_\alpha({\bf x})=x_\alpha, g_\beta({\bf x})=1-x_\beta$,
inequality \apBKR of Theorem \ref{thm1-dual} gives \bea
\label{X(1-X)} E \left\{ \max_{\alpha \not = \beta}X_\alpha
(1-Y_\beta) \right\} \le E X_{[1]}(1-X_{[n]}).\ena Note that
$\max_{\alpha \not = \beta}X_\alpha (1-Y_\beta) \ne X_{[1]}
(1-Y_{[n]})$; the right hand side might be larger because of the
restriction $\alpha \not =\beta$. Removing the restriction $\alpha
\not = \beta$ reverses (\ref{X(1-X)}), that is,
$$ E X_{[1]}(1-X_{[n]}) \le E X_{[1]} E(1-X_{[n]}) = E X_{[1]}
E(1-Y_{[n]})=E X_{[1]} (1-Y_{[n]}) =
E\left\{\max_{\alpha,\beta}X_\alpha (1-Y_\beta)\right\},$$ where
the inequality follows by the negative association of $X_{[1]}$
and $1-X_{[n]}$.

Following our treatment of applications of Theorem \ref{thm1} we
can extend (\ref{X(1-X)}) as follows: with ${\cal A}$ and ${\cal
B}$ the collection of all $k$ and $l$ subsets of $\N$
respectively, and
$$ f_\alpha({\bf x})=\prod_{j \in \alpha}x_j \quad \mbox{and}
\quad g_\beta({\bf x})=\prod_{j \in \beta}(1-x_j), $$ we obtain
$$ E \left\{ \max_{\alpha
\cap \beta =\emptyset}\prod_{i \in \alpha, j \in \beta}X_i (1-Y_j)
\right\} \le E \left\{ \prod_{1 \le i \le k, 1 \le j \le l}
X_{[i]} (1-X_{[n-j+1]}) \right\}. $$
\end{example}

We now consider resource allocation problems of the following
type. Suppose that two projects $A$ and $B$ have to be completed
using $n$ available resources represented by the components of a
vector ${\bf x}$. Each resource can be used for at most one
project, and an allocation is given by a specification of disjoint
subsets of resources. For any given subsets $\alpha,\beta
\subseteq \N$, let $f_\alpha({\bf x})$ and $g_\beta({\bf x})$
count the number of ways that projects $A$ and $B$ can be
completed using the resources ${\bf x}_\alpha$ and ${\bf x}_\beta$
respectively. The exact definitions of the projects and the counts
are immaterial; in particular larger sets do not necessarily imply
more ways to carry out a project. For an  allocation
$\alpha,\beta$, $\alpha \cap \beta= \emptyset$,  the total number
of ways to carry out the two projects together is the product
$f_\alpha({\bf x})g_\beta({\bf x})$. When the resources are
independent variables, inequality \aBKR\, of Theorem \ref{thm1}
bounds the expected maximal number of ways of completing $A$ and
$B$ together, by the product of the expectations of the maximal
number of ways of completing each project alone. The bound is
simple in the sense that it does not require understanding of the
relation between the two projects. In particular, it can be
computed without knowledge of the optimal allocation of resources.

\begin{example}
With $J$ a list of tasks, consider fulfilling the set of tasks on
(not necessarily disjoint) lists $A \subseteq J$ and $B \subseteq
J$, in two distant cities using disjoint sets of workers chosen
from $1,2,\ldots,n$. Each worker may be sent to one of the cities
and assigned a single task or a set of tasks which he can perform.
A worker may be qualified to fill more than one set of tasks. For
$i \in \N$, let ${\bf x}_i \subseteq 2^J$ be the collection of
possible assignments of tasks for worker $i$, (that is, the sets
of tasks worker $i$ is qualified to fulfill.); For $\alpha, \beta
\subseteq \N$ and ${\bf x}=({\bf x}_1,\ldots,{\bf x}_n)$, let
$f_\alpha({\bf x})$ equal the number of ways the collection of
workers $\alpha$ can complete $A$, and $g_\beta({\bf x})$ the
number of ways the collection $\beta$ can complete $B$. When the
qualifications ${\bf X}_i, i \in \N$ are independent, Theorem
\ref{thm1} bounds the expectation of the maximal number of ways of
fulfilling the task requirements in both cities, by the product of
the expectations of the maximal numbers of ways that the
requirements in each collection can be separately satisfied.
\end{example}

\begin{example}
\textbf{Paths on Graphs} Consider a graph $\cal G$ with an
arbitrary fixed vertex set ${\cal V}=\{v_1,\ldots,v_n\}$, where
for each pair of vertices the existence of the edge $\{v_i,v_j\}$
is determined independently using a probability rule based on
$v_i,v_j$, perhaps depending only on $d(\{v_i,v_j\})$ for some
function $d$. Let ${\bf X}=\{X_{\{i,j\}}\}$ where $X_{\{i,j\}}$ is
the indicator that there exists an edge between $v_i$ and $v_j$.
For instance, with ${\cal V} \subseteq {\bf R}^m$  and
$Z_{\{i,j\}}, 1 \le i, j \le n$ independent non-negative
variables, we may take for ${ v}_i, {v}_j \in {\cal V}$, \beas
X_{\{i,j\}}={\bf 1}(d(\{{ v}_i,{ v}_j\})<Z_{\{i,j\}}) \quad
\mbox{where} \quad d(\{{ v}_i,{ v}_j\})=||{ v}_i-{ v}_j|| \quad
\mbox{(Euclidean distance)} \quad.
\enas Note that since the variables $Z_{\{i,j\}}$ do not have to be
identically distributed, we can set $Z_{i,i}=0$ and avoid self loops
should we wish to do so.

Let a path in the graph $\cal G$ from $u$ to $w$ be any ordered
tuple of vertices $v_{i_1},\ldots,v_{i_p}$ with $v_{i_1}=u,
v_{i_p}=w$ and $X_{\{i_k,i_{k+1}\}}=1$ for $k=1,\ldots,p-1$, and
having all edges $\{v_{i_k},v_{i_{k+1}}\}$ distinct. For $u,v$ and
$w$ in $\cal V$ and $\alpha, \beta \subseteq \{\{i,j\}: 1 \le i, j
\le n \}$, let $f_\alpha({\bf X})$ be the number of paths in the
graph from $u$ to $v$ which use only edges $\{v_i,v_j\}$ for
$\{i,j\} \in \alpha$; in the same manner, let $g_\beta({\bf X})$ be
the number of paths in the graph from $v$ to $w$ which use only
edges $\{v_i,v_j\}$ for $\{i,j\} \in \beta$.

 The ``projects" $A$ and $B$ in this framework are to create
 paths from $u$ to $v$ using $\alpha$,
and from $v$ to $w$ using $\beta$, respectively, which combine
together, when $\alpha \cap \beta = \emptyset$, to give the overall
project of creating a path from $u$ to $w$ passing through $v$. As
the product $f_\alpha({\bf X})g_\beta({\bf X})$ for $\alpha \cap
\beta = \emptyset$ is the number of paths from $u$ to $w$ via $v$
for the given allocation, Theorem \ref{thm1} provides a bound on the
expected maximal number of such paths over all allocations in terms
of the product of the expectations of the maximal number of paths
from $u$ to $w$ and from $w$ to $v$ when the paths are created
separately. Though finding the optimal allocation may be demanding,
the upper bound can be computed simply, for this case in particular
by monotonicity of $f_\alpha({\bf x}), g_\beta({\bf x})$ in $\alpha$
and $\beta$ for fixed ${\bf x}$, implying that the maximal number of
paths created separately is attained when using all available edges,
i.e. at $\alpha = \beta =\N$.

However, the result and the upper bound hold even in constrained
situations where the existence of more edges does not lead to more
paths, that is, in cases where the functions $f_\alpha, g_\beta$ are
not monotone in $\alpha$ and $\beta$. One such case would be where
the existence of a particular edge mandates that all paths from $u$
to $v$ use it. More specifically, for some fixed $\{i_0,j_0 \}$
suppose that if $\{i_0,j_0 \} \in \alpha$ and $x_{\{i_0,j_0\}}=0$
then $f_\alpha({\bf x})$ counts the number of paths from $u$ to $v$.
On the other hand if $x_{\{i_0,j_0\}}=1$ then $f_\alpha({\bf x})$
counts the number of paths from  $u$ to $v$ which use the edge
$\{v_{i_0},v_{j_0}\}$. In general such $f_\alpha$ will not be
monotone.

This example easily generalizes to paths with multiple waypoints.
We may also consider directed graphs where for $1 \le i \not =j
\le n$ the directed edge $(v_i,v_j)$ from $v_i$ to $v_j$ exists
when $X_{ij}=1$, the directed edge $(v_j,v_i)$ from $v_j$ to $v_i$
exists when $X_{ij}=-1$ and $X_{ij}=0$ when no edge exists.
Returning to the graph example following the statement of Theorem
\ref{BKR}, when the signed edge indicators $\{X_{ij}\}_{1 \le i<j
\le n}$ are independent, inequality \aBKR\, of Theorem \ref{thm1}
provides a bound on the expected maximal number of paths from
vertices $v_1$ to $v_2$ and $w_1$ to $w_2$ using disjoint edges.
Another possible extension is to consider paths between subsets of
vertices.

For application of the dual inequality, consider for example two
directed graphs on the same vertex set, determined by equally
distributed and independent collections of signed edge indicators
$\bf X$ and $\bf Y$, each having independent (but not necessarily
identically distributed) components. Let $\alpha, \beta \subseteq
\{(i,j): 1 \le i \not = j \le n \}$, and $f_\alpha({\bf X})$ be
the number of directed paths in the graph from vertices $u$ to $v$
which use only $\bf X$ edges $(v_i,v_j)$ with $(i,j) \in \alpha$;
in the same manner, let $g_\beta({\bf Y})$ be the number of
directed paths in the graph from $v$ back to $u$ which use only
$\bf Y$ edges $(v_i,v_j)$ with $(i,j) \in \beta$. Consider the
expected maximal number of paths, over all $\alpha$ and $\beta$
with $\alpha \cap \beta = \emptyset$, that go from $u$ to $v$
using the ${\bf X}$ edges $\alpha$ and return to $u$ from $v$
using the ${\bf Y}$ edges $\beta$. Then Theorem \ref{thm1-dual}
implies that this expectation is bounded by the expected maximal
number of paths, over all $\alpha$ and $\beta$, to move from $u$
to $v$ using $\alpha$, and then returning to $u$ using $\beta$,
all with ${\bf X}$ edges, but where edges used on the forward trip
may now also be used for the return.
\end{example}

\section{Proofs}
\label{proof}
\subsection{Proofs of Proposition \ref{prop1} and Theorem \ref{thm1}}
We first reduce the problem by proving the following implications
between the parts of Theorems \ref{thm1} and Proposition
\ref{prop1}.
\begin{proposition}
\label{reformulations} $\aBKR \Rightarrow \cBKR \Rightarrow \dBKR
\Rightarrow \aBKR$ and $\bBKR \Leftrightarrow \eBKR$.
\end{proposition}

\noindent {\bf Proof:} $\aBKR \Rightarrow \cBKR$: Apply inequality
\aBKR\, to the finite collections $\{\tilde{f}_K\}_{K \in {\cal
K}},\{\tilde{g}_L\}_{L \in {\cal L}}$ and use \bea \nonumber
\sup_{\alpha \cap \beta = \emptyset}f_\alpha({\bf x})g_\beta({\bf
x})=\max_{\stackrel{K \cap L = \emptyset}{K \in {\cal K}, L \in
{\cal L}}} \left(\sup_{K_\alpha \subseteq K, L_\beta
\subseteq L}f_\alpha({\bf x})g_\beta({\bf x}) \right)\\
\label{I-to-II} = \max_{\stackrel{K \cap L = \emptyset}{K \in
{\cal K}, L \in {\cal L}}} \left( \sup_{K_\alpha \subseteq K}
f_\alpha({\bf x})\sup_{L_\beta \subseteq L}g_\beta({\bf
x})\right)=\max_{\stackrel{K \cap L = \emptyset}{K \in {\cal K}, L
\in {\cal L}}} {\tilde f}_K({\bf x}){\tilde g}_L({\bf x}). \ena

\noindent $\cBKR \Rightarrow  \dBKR$: Apply $\tilde{f} _{K({\bf
x})}({\bf x}) \tilde{g}_{L({\bf x})}({\bf x}) \le \max_{K \cap L =
\emptyset, K \in {\cal K}, L \in {\cal L}}\tilde{f} _{K}({\bf x}) \tilde{g}_{L}({\bf x}).$\\[0.5ex]
\noindent $\dBKR \Rightarrow  \aBKR$: Note that the right hand
side of (\ref{I-to-II}) equals $\tilde{f}_{K({\bf x})}({\bf x})
\tilde{g}_{L({\bf x})}({\bf x})$ for some $K({\bf x}) \in {\cal
K}$ and $L({\bf
x})\in {\cal L}$ with $K({\bf x}) \cap L({\bf x}) = \emptyset$.\\[0.5ex]
\noindent $\bBKR \Rightarrow \eBKR$: Apply $\underline{f} _{K({\bf
x})}({\bf x}) \underline{g}_{L({\bf x})}({\bf x}) \le \max_{K \cap
L = \emptyset, K \in {\cal K}, L \in {\cal L}}\underline{f}
_{K}({\bf x}) \underline{g}_{L}({\bf
x}).$\\[0.5ex]
\noindent $\eBKR \Rightarrow \bBKR$: Use the fact that there exist
some disjoint $K({\bf x}) \in {\cal K},L({\bf x}) \in {\cal L}$
such that \bea \label{maxKL-measurable} \max_{\stackrel{K \cap L =
\emptyset}{K \in {\cal K}, L \in {\cal L}}} \underline{f}_K({\bf
x})\underline{g}_{L}(\mathbf x) = \underline{f}_{K(\bf x)}({\bf
x})\underline{g}_{L({\mathbf x})}(\mathbf x). \quad \qed \ena

Let ${\cal F}_C$ be the sigma algebra generated by a collection of sets $C$. We say ${\cal F}_C$ is a finite product \lcolor{sigma sub algebra of $\mathbb{S}$ when}
\bea \label{eq:calC}
\lcolor{{\cal C}=\left\{\prod_{i=1}^n A_i, A_i \in {\cal C}_i\right\},
\qmq{with
${\cal C}_i\subseteq \mathbb{S}_i$ finite for all
$i=1,\ldots,n.$}}
\ena
\lcolor{It is easy to see that every finite sigma algebra, ${\cal F}$, contains a subset $G$, not containing the empty set, such that every element of ${\cal F}$ can be represented uniquely as a disjoint union of elements of $G$. Call $G$ the disjoint generating set of ${\cal F}$.}

Our next objective is proving the inequalities of Framework \ref{I}, to be accomplished by proving \dBKR\, in Lemma \ref{finite=infinite}.
We start with a simple extension of inequality (\ref{original-box}),
expressed in terms of indicator functions, \lcolor{from finite spaces to spaces that may not be finite, but which are endowed with a finite product sigma algebra. }

\begin{lemma}
\label{almost-trivial} Let $Q$ be any probability product measure
on the finite product sigma algebra ${\cal F}_{\cal C}$ with ${\cal C}$ given by \eqref{eq:calC}. Then, inequality
\aBKR\, holds when expectations are taken with respect to $Q$, and
$\{f_\alpha\}_{\alpha \in {\cal A}},\{g_\beta\}_{\beta \in {\cal
B}}$ are ${\cal F}_C$ measurable indicator functions.
\end{lemma}

\noindent {\bf Proof:} \lcolor{For $i=1,\ldots,n$, let $G_i$ be the disjoint generating set of ${\cal F}_{{\cal C}_i}$. By Theorem \ref{BKR}, applied on the space $G=\prod_{i=1}^n G_i$,
\bea
\label{QAB}
Q(A \Box B) \le Q(A)Q(B).
\ena }
Let events $A$ and $B$ be defined by the indicator functions
\bea
\label{not-5}
{\bf 1}_A({\bf x})=\max_\alpha f_\alpha({\bf x}),
\quad {\bf 1}_B({\bf x})=\max_\beta g_\beta({\bf x}),
\ena
\lcolor{and let $A_\alpha$ and $B_\beta$ be the sets indicated by $f_\alpha({\bf x})$ and $g_\beta({\bf x})$ respectively. Suppose ${\bf x} \in S$ satisfies $f_\alpha({\bf x})g_\beta({\bf x})=1$ for disjoint $\alpha,\beta$. Clearly $A_\alpha \subseteq A$, and as $f_\alpha$ depends on $K_\alpha$, we have $[{\bf x}]_\alpha \subseteq A_\alpha \subseteq A$. As a similar statement holds for $B$, ${\bf x} \in A \Box B$, hence, }
\bea \label{boximpliesbkrforindicators} \max_{ \alpha \cap
\beta = \emptyset} f_\alpha({\bf x})g_\beta({\bf x}) \leq{\bf 1}_{A \Box B}({\bf x}). \ena

Now (\ref{boximpliesbkrforindicators}) gives the
first inequality below, (\ref{QAB}) the second inequality, and
(\ref{not-5}) the last equality in \beas E_Q \left\{\max_{ \alpha
\cap \beta = \emptyset} f_\alpha({\bf X})g_\beta({\bf X})\right\}
\le Q(A \Box B) \le Q(A)Q(B)= E_Q \left\{\max_\alpha f_\alpha({\bf
X})\right\} \,E_Q\left\{\max_\beta g_\beta({\bf X})\right\}. \qed
\enas

We say a collection of functions is FP if it generates a finite product sigma algebra contained in $\mathbb S$; note that a finite union of FP collections is FP.

\begin{lemma}
\label{base-step-S-null} Inequality \aBKR\, is true for $P$ any
probability  product measure on $({\spa}, {\mathbb S})$, and
$\{f_\alpha\}_{\alpha \in {\cal A}},\{g_\beta\}_{\beta \in {\cal
B}}$, any finite collections of FP indicator functions.
\end{lemma}

\noindent {\bf Proof:} Let ${\cal H}$ be the sigma algebra
generated by $\{f_\alpha\}_{\alpha \in {\cal
A}},\{g_\beta\}_{\beta \in {\cal B}}$, and $Q :=P|_{\cal H}$, the
restriction of $P$ to the finite product sigma algebra ${\cal H}$.
For $h$ an ${\cal H}$ measurable indicator function, that is, for
$h({\bf x})={\bf 1}_A({\bf x})$ for some $A \in {\cal H}$, we have
\bea \label{same-measure} E_Q h = Q(A) = P(A)=E_P h. \ena Since
the product of ${\cal H}$ measurable indicators is an ${\cal H}$
measurable indicator, and the same is true for the maximum, we
have by Lemma \ref{almost-trivial} and (\ref{same-measure}), \beas
E_P \left\{\max_{ \alpha \cap \beta = \emptyset} f_\alpha({\bf
X})g_\beta({\bf X})\right\} &=& E_Q \left\{\max_{ \alpha \cap
\beta = \emptyset} f_\alpha({\bf X})g_\beta({\bf X})\right\}\\
&\le& E_Q \left\{\max_\alpha f_\alpha({\bf X})\right\}
\,E_Q\left\{\max_\beta g_\beta({\bf X})\right\}\\ &=&
E_P\left\{\max_\alpha f_\alpha({\bf X})\right\} \,E_P
\left\{\max_\beta g_\beta({\bf X})\right\}. \,\,\quad \qed\enas

Let ${\cal P}$ denote the collection of all product sets of the
form ${\cal C}=\{\prod_{i=1}^n S_i, S_i \in {\cal C}_i\}$ where
${\cal C}_i\subseteq \mathbb{S}_i$ are finite for all
$i=1,\ldots,n$. Then
\bea
\label{Jay} {\mathbb S}=\cal F_{\cal J}  \quad \mbox{where} \quad
{\cal J}=\bigcup_{{\cal C} \in {\cal P}}{\cal F}_{\cal C}. \ena

Lemma \ref{inductive-step-S-null} generalizes the inequality from FP
indicator functions to ${\mathbb S}$ measurable indicator functions.
\begin{lemma}
\label{inductive-step-S-null} Inequality \aBKR\, is true for any
probability product measure $\PP$ and finite collections
$\{f_\alpha\}_{\alpha \in {\cal A}}$ and $\{g_\beta\}_{\beta \in
{\cal B}}$ of\, ${\mathbb S}$ measurable indicator functions.
\end{lemma}
\noindent {\bf Proof}:
For ${\cal R},{\cal S}$ subsets of ${\cal A} \cup {\cal B}$
satisfying ${\cal R} \cap {\cal S}=\emptyset$ and ${\cal R} \cup
{\cal S}={\cal A} \cup {\cal B}$, we proceed by induction on the
cardinality of the set ${\cal S}$ in the statement $I({\cal
R},{\cal S})$: inequality \aBKR\, is true when $\{f_\alpha\}_{\alpha
\in {\cal R}\cap {\cal A}} ,\{g_\beta\}_{\beta \in {\cal R} \cap
{\cal B}}$ are finite FP collections of indicator functions, and
$\{f_\alpha\}_{\alpha \in {\cal S} \cap {\cal
A}},\{g_\beta\}_{\beta \in {\cal S} \cap {\cal B}}$ are any finite
collections of ${\mathbb S}$ measurable indicators. Lemma
\ref{base-step-S-null} shows that $I({\cal A} \cup {\cal
B},\emptyset)$ is true, and the conclusion of the present lemma is
$I(\emptyset,{\cal A} \cup {\cal B})$. Assume for some such ${\cal
R},{\cal S}$ with ${\cal S} \not = {\cal A} \cup {\cal B}$, that
$I({\cal R},{\cal S})$ is true. For $\gamma \in {\cal R}$ with,
say $\gamma \in {\cal A}$, let ${\cal M}$ be the collection of all
sets $A \subseteq \spa$ such that \aBKR\, holds for $f_\gamma={\bf
1}_A$, and when $\{f_\alpha\}_{\alpha \in {\cal R} \cap {\cal A}
\setminus \{\gamma\} }$ and $\{g_\beta\}_{\beta \in {\cal R} \cap {\cal
B}}$ are finite FP indicators, and
$\{f_\alpha\}_{\alpha \in {\cal S} \cap {\cal
A}},\{g_\beta\}_{\beta \in {\cal S}\cap {\cal B}}$ are any collection
of ${\mathbb S}$ measurable indicators. The singleton collection
$f_\gamma$ is FP for any $A \in {\cal J}$ given by (\ref{Jay}).
Therefore, for any $A \in {\cal J}$, the union $f_\gamma,
\{f_\alpha\}_{\alpha \in {\cal R} \cap {\cal A} \setminus
\{\gamma\} },\{g_\beta\}_{\beta \in {\cal R}\cap {\cal B}}$ is FP.
By the induction hypothesis, ${\cal J} \subseteq {\cal M}$. Since
${\cal M}$ is a monotone class and ${\cal J}$ is an algebra which
generates ${\mathbb S}$, the monotone class theorem implies
${\mathbb S} \subseteq {\cal M}$. This completes the induction.
$\qed$

We now relax the requirement that the functions be indicators.
\begin{lemma}
\label{simple} Inequality \dBKR \,\,is true for any product
measure $\PP$ and finite collections of ${\mathbb S}$ measurable
functions $\{f_\alpha\}_{\alpha \in {\cal A}}$ and
$\{g_\beta\}_{\beta \in {\cal B}}$ which assume finitely many
non-negative values.
\end{lemma}
\noindent {\bf Proof}: We prove \dBKR \,\,by induction on $m$ and $l$, the
number of values taken on by the collections $\{f_\alpha\}_{\alpha
\in {\cal A}}, \{g_\beta\}_{\beta \in {\cal B}}$, respectively. By
Lemma \ref{inductive-step-S-null} inequality \aBKR\, is true for
finite collections of measurable indicators, and hence by
Proposition \ref{prop1}, so is inequality \dBKR. Now the base case
$m=2$, $l=2$ follows readily by extending from indicators to two
valued functions by linear transformation.

Assume the result is true for some $m$ and $l$ at least 2, and
consider a collection $\{f_\alpha\}_{\alpha \in {\cal A}}$
assuming the values $0 \le a_1<\cdots<a_{m+1}$; a similar argument
applies to induct on $l$. For some $k$, $2 \le k \le m$, define
$$
A_{\alpha,k}=\{{\bf x}: f_\alpha({\bf x})=a_k\},
$$
and for $a_{k-1} \le a \le a_{k+1}$, let
$$
h_\alpha^a({\bf x})=f_\alpha({\bf x})+(a-a_k){\bf
1}_{A_{\alpha,k}}({\bf x}),
$$
the function $f_\alpha$ with the value of $a_k$ replaced by $a$.
We shall prove that for all $a \in [a_{k-1},a_{k+1}]$ inequality
\dBKR \,\, holds with $\{h_\alpha^a\}_{\alpha \in {\cal A}}$
replacing $\{f_\alpha\}_{\alpha \in {\cal A}}$. By the induction
hypothesis we know it holds at the endpoints, that is, for $a \in
\{a_{k-1},a_{k+1}\}$, since then the collection
$\{h_\alpha^a\}_{\alpha \in {\cal A}}$ takes on $m$ values;
clearly, the case $a=a_k$ suffices to prove the lemma.

Given $\Gamma({\bf x})$, a function with values in $2^{\cal A}$,
with some abuse of notation denote \bea \label{abuse}
\tilde{f}_{\Gamma({\bf x})}({\bf x}) = \sup_{\alpha: \alpha \in
\Gamma({\bf x})} f_\alpha({\bf x}). \ena Note that
$\tilde{f}_{K({\bf x})}({\bf x})$ in (\ref{Kalphax}) corresponds
to $\Gamma({\bf x})=\{\alpha: K_\alpha \subseteq K({\bf x}) \}$,
and similarly for $\tilde{g}_{L({\bf x})}({\bf x})$; for
measurability issues see Section \ref{appendix}. For any
function $\Gamma({\bf x})$ with values in $2^{\cal A}$, we have
for all $a \in [a_{k-1},a_{k+1}]$, \beas
C_\Gamma:=\{ {\bf x}: \widetilde{h^a}_{\Gamma({\bf x})} ({\bf
x})=a, \, \tilde{f}_{\Gamma({\bf x})}({\bf x}) \not \in
\{a_{k-1},a_{k+1}\}\} = \{{\bf x}: \tilde{f}_{\Gamma({\bf
x})}({\bf x})=a_k\}, \enas showing that $C_\Gamma$ does  not
depend on $a$.

Let $D=C_\Gamma$ for $\Gamma({\bf x})={\cal A}$, and note that
$$ \sup_\alpha h^a_\alpha({\bf x}) = a 1_D + \sup_\alpha
f_\alpha({\bf x}) 1_{D^c}.$$ Then the right hand side of \dBKR,
with $\{h_\alpha^a\}_{\alpha \in {\cal A}}$ replacing
$\{f_\alpha\}_{\alpha \in {\cal A}}$, equals $a\delta+\lambda$,
where
$$
\delta=\PP(D)\int \sup_\beta g_\beta(\mathbf x) d \PP(\mathbf x)
\quad \mbox{and} \quad \lambda=\int_{D^c}\sup_\alpha f_\alpha({\bf
x})d \PP(\mathbf x)\int \sup_\beta g_\beta(\mathbf x) d \PP(\mathbf
x)
$$
do not depend on $a$. Now, let $E=C_\Gamma$ for $\Gamma({\bf
x})=\{\alpha: K_\alpha=K({\bf x})\}$ and note that
$\widetilde{h^a}_{{ K}({\bf x})}({\bf x})=a1_E + \widetilde{f}_{{
K}({\bf x})}({\bf x})1_{E^c}$. Similarly, the left hand side of
\dBKR, with $\{h_\alpha^a\}_{\alpha \in {\cal A}}$ replacing
$\{f_\alpha\}_{\alpha \in {\cal A}}$, equals $a\theta+\eta$, where
$$
\theta=\int_E \tilde{g}_{L({\mathbf x})} (\mathbf x)d\PP(\mathbf x)
\quad \mbox{and} \quad \eta= \int_{E^c}\widetilde{f}_{{ K}({\bf
x})}({\bf x})\tilde{g}_{{ L}({\bf x})}({\bf x})d\PP(\mathbf x)
$$
do not depend on $a$.  When $a \in \{a_{k-1},a_{k+1}\}$ the
collection $h_\alpha,\alpha \in {\cal A}$ takes on $m$ values, so
by the induction hypotheses \dBKR\, holds with
$\{h_\alpha^a\}_{\alpha \in {\cal A}}$ replacing
$\{f_\alpha\}_{\alpha \in {\cal A}}$ and we obtain \bea
\label{take-cc} a \theta + \eta \le a \delta + \lambda, \quad
\mbox{for $a \in \{a_{k-1},a_{k+1}\}$.} \ena By taking a convex
combination, we see that inequality (\ref{take-cc}) holds for all
$a \in [a_{k-1},a_{k+1}]$, so in particular for $a_k$, completing
the induction. $\qed$

\begin{lemma}
\label{simpler} Inequality \dBKR \,\, is true for any probability
product measure $\PP$ and finite collections of non-negative
${\mathbb S}$ measurable functions $\{f_\alpha\}_{\alpha \in {\cal
A}}$ and $\{g_\beta\}_ {\beta \in {\cal B}}$.
\end{lemma}

\noindent {\bf Proof:} Lemma \ref{simple} shows that the result is
true for simple functions. By approximating the functions
$f_\alpha,g_\beta$ below by simple functions, $f_{\alpha,k}
\uparrow f_\alpha, g_{\beta,k} \uparrow g_\beta$ as $k \uparrow
\infty$, and applying the monotone convergence theorem, we have
the result for arbitrary non-negative functions.$\qed$

\begin{lemma}
\label{finite=infinite} Inequality \dBKR \, is true for countable
collections of non-negative ${\mathbb S}$ measurable functions
$\{f_\alpha\}_{\alpha \in {\cal A}}$ and $\{g_\beta\}_ {\beta \in
{\cal B}}$.
\end{lemma}

\noindent {\bf Proof:} For $K, L \in 2^{\bf n}$, let
$$
\varphi_K({\bf x})=\tilde{f}_K({\bf x}) \quad \mbox{and} \quad
\phi_L({\bf x})=\tilde{g}_K({\bf x}),
$$
recalling definition (\ref{one-to-2}). Noting
$$
\tilde{f}_{K(\bf x)}({\bf x})=\sup_{\alpha:K_\alpha \subseteq
K({\bf x})}f_\alpha({\bf x})=\sup_{K \subseteq K({\bf
x})}\sup_{\alpha:K_\alpha \subseteq K}f_\alpha({\bf x})=\sup_{K
\subseteq K({\bf x})}\tilde{f}_K({\bf x})=\tilde{\varphi}_{K({\bf
x})}({\bf x}),
$$
and
\beas  \sup_K \varphi_K({\bf x})=\sup_K \tilde{f}_K({\bf x})
= \sup_\alpha f_\alpha({\bf x}), \enas and similarly for
$\{g_\beta\}_{\beta \in \cal B}$, the result follows immediately
upon applying Lemma \ref{simpler} to the finite collections
$\{\varphi_K\}_{K \in 2^\N}$ and $\{\phi_L\}_{L \in 2^\N}$. $\qed$

By Proposition \ref{reformulations}, at this point we have completed
proving all inequalities pertaining to Framework \ref{I}. The next
proposition connects the two frameworks and completes the proof of
Theorem \ref{thm1}, and again applying Proposition
\ref{reformulations}, that of Proposition \ref{prop1}.
\begin{proposition}
Inequality \aBKR\, holds in Framework \ref{I} for all
collections $\{f_\alpha\}_{\alpha \in {\cal A}},\{g_\beta\}_{\beta
\in {\cal B}}$ of given functions, if and only if inequality \bBKR\,
holds in Framework \ref{II} for all given functions $f$ and $g$ and
collections $\cal K$ and $\cal L$.
\end{proposition}
\noindent {\bf Proof:} $\aBKR \Rightarrow \bBKR$. \lcolor{For $L \subseteq \N$ let $P_L({\bf x})$ denote the marginal of $P$ in the coordinates indexed by $L$.
Let functions $f$, $g$ and collections $\cal K$ and $\cal L$ of subsets of $2^{\bf n}$
be given.}

\lcolor{By Fubini's theorem, for any $K \subseteq 2^{\bf n}$,
\begin{multline*}
P(\underline{f}_K({\bf x}) \le f({\bf x})) = \int {\bf 1}(\underline{f}_K({\bf x}) \le f({\bf x}))dP({\bf x})
=\int \int  {\bf 1}(\underline{f}_K({\bf x}) \le f({\bf x}))dP_{K^c}({\bf x})dP_K({\bf x}) \\= \int P_{K^c}(\underline{f}_K({\bf x}) \le f({\bf x})) dP_K({\bf x})= \int 1 dP_K({\bf x})=1,
\end{multline*}
where the fourth equality holds by definition of the essential infimum.
 As ${\cal K}$ is finite,
	\beas
	P(\max_{K \in {\cal K}} \underline{f}_K({\bf X}) \le f({\bf X}))=1,
	\quad \mbox{implying} \quad
E\left\{
	\max_{K \in {\cal K}} \underline{f}_K({\bf X} \right\}
	\le E\left\{ f({\bf X}) \right\},
\enas
with a similar inequality holding for $g$. Now we see that \bBKR\, holds by applying \aBKR\, to the
}
collections
$\{\underline{f}_K({\bf x})\}_{K \in \cal K}$ and
$\{\underline{g}_L({\bf x})\}_{L \in \cal L}$ as in
(\ref{fgismin}).

{\noindent $\bBKR \Rightarrow  \aBKR$: Given  collections
of functions $f_\alpha,g_\beta $ depending on $K_\alpha,L_\beta$,
define
\begin{equation}
\label{ssup} f({\bf x}) = \sup_\alpha f_{\alpha}({\bf x}) \quad
\mbox{and} \quad g({\bf x}) = \sup_\beta g_{\beta}({\bf x}).
\end{equation}

Now letting $\underline{f}_K, \underline{g}_L$ be as in
(\ref{fgismin}), we have \beas 
f_{\alpha}({\bf x}) = \underline{f_{\alpha}}_{K_\alpha}({\bf x})
\le \underline{f}_{K_\alpha}({\bf x}) \quad \mbox{and likewise}
\quad g_{\beta}({\bf x})\le \underline{g}_{L_\beta}({\bf x}).
\enas

Now, for $\alpha, \beta$ disjoint,
\beas
f_{\alpha}({\bf x})g_{\beta}({\mathbf x}) \le
\underline{f}_{K_\alpha}({\bf x})\underline{g}_{L_\beta}({\mathbf
x})  \le \max_{\stackrel{K \cap L = \emptyset}{K \in 2^{\bf n}, L
\in 2^{\bf n}}} \underline{f}_K({\bf x})\underline{g}_{L}(\mathbf
x). \enas Taking supremum on the left hand side over all disjoint
$\alpha,\beta$ and then expectation, the result now follows by
applying inequality \bBKR\, and (\ref{ssup}).}\,\,\quad \qed

\subsection{The Dual Inequality}

As observed in [\ref{kss}], the techniques in [\ref{vbf}] extend
the dual inequality (\ref{dual-reimer}) from uniform measure on
$\{0,1\}^n \times \{0,1\}^n$ to any product measure on a discrete
finite product space $\spa$. Specifically, Lemmas 3.2(iii), 3.4,
and 3.5 of [\ref{vbf}] carry over with minimal changes,
essentially by replacing $\Box$ by $\Diamond$ and $\cap$ by
$\times$ appropriately; for example, the dual version of Lemma 3.4
would begin with the identity
\beas (f \times f)^{-1}(A \Diamond B) = \bigcup_{C_1,C_2}
\left\{(f \times f)^{-1}(C_1 \times C_2)\right\} \enas where the
union is over all $C_1,C_2$ such that $C_1$ is a maximal cylinder
of $A$, $C_2$ is a maximal cylinder of $B$, and $C_1 \perp C_2$;
see Sections 3 and 2 of [\ref{vbf}] for the formal definitions of
maximal cylinder, and perpendicularity $\perp$, respectively.

Now the proof of Theorem \ref{thm1-dual} and Proposition
\ref{prop1-dual} follow in a nearly identical manner to that of
Theorem \ref{thm1} and Proposition \ref{prop1}.  For instance, to
prove \dpBKR, consider \beas
C_\Gamma=\{ ({\bf x},{\bf y}): \widetilde{h^a}_{\Gamma({\bf
x},{\bf y})} ({\bf x})=a, \, \tilde{f}_{\Gamma({\bf x},{\bf
y})}({\bf x}) \not \in \{a_{k-1},a_{k+1}\}\} = \{({\bf x},{\bf
y}): \tilde{f}_{\Gamma({\bf x},{\bf y})}({\bf x})=a_k\}. \enas
Setting $D=C_\Gamma$ for $ \Gamma({\bf x},{\bf y})=\{\alpha:
K_\alpha=K({\bf x},{\bf y})\} $ we can write the left hand side of
\dpBKR as $a\theta+\eta$, with
$$
\theta= \int_D \tilde{g}_{L({\bf x},{\bf y})}({\bf y})d\PP({\bf
x})d\PP({\bf y}) \quad \mbox{and} \quad
\eta=\int_{D^c}\widetilde{f}_{K({\bf x},{\bf y})}({\bf
x})\tilde{g}_{L({\bf x},{\bf y})}({\bf y})d\PP({\bf x})d\PP({\bf
y}),
$$
and using $E=C_\Gamma$ for $\Gamma={\cal A}$, the right hand side
becomes $a\delta+\lambda$, where
$$
\delta = \int_D \sup_\beta g_\beta({\bf x})d\PP({\bf x}) \quad
\mbox{and} \quad \lambda =
\int_{D^c}\sup_{\alpha,\beta}f_\alpha({\bf x}) g_\beta({\bf
x})d\PP({\bf x})
$$
with $\theta,\eta,\delta$ and $\lambda$ not depending on $a$.

\section{A PQD ordering inequality}
\label{PQD-section} Consider a collection $\{f_\alpha({\bf
x})\}_{\alpha=1}^m$ of functions which are all increasing or all
decreasing in each component of ${\bf x}=(x_1, \ldots, x_n) \in
{\bf R}^n$. Let ${\bf X}=(X_1,\ldots,X_n) \in {\bf R}^n$ be a
vector of independent random variables, ${\bf Y}=(Y_1,\ldots,Y_n)$
an independent copy of ${\bf X}$, and for each
$\alpha=1,\ldots,m$, let $H_\alpha \subseteq \N$, and \bea
\label{defYj} {\bf Z}_\alpha=(Z_{1,\alpha},\ldots,Z_{n,\alpha}),
\ena where $Z_{i,\alpha} = Y_i$ if  $i \in H_\alpha$, and
$Z_{i,\alpha}= X_i$, if $i \not \in H_\alpha$. Now let \bea
\label{defSk} \vecS=(f_1({\bf Z}_1),\ldots,f_m({\bf Z}_m))\quad
\mbox{and} \quad \vecR=(f_1({\bf X}),\ldots,f_m({\bf X})). \ena
Inequalities between vectors below are coordinate-wise. When
(\ref{PQDinq}) below holds, we say that the components of $\vecR$
are more `Positively Quadrant Dependent' than those of $\vecS$,
and write $\vecS \le_{PQD} \vecR$.

\begin{theorem}
\label{mon} For every ${\bf c}=(c_1,\ldots,c_m) \in {\bf R}^m$ and
$H_\alpha \subseteq {\bf n}, \alpha=1,\ldots,m$,
\bea \label{PQDinq}
  P(\vecS \ge {\bf c}) \le P(\vecR \ge {\bf c}) \quad
  \mbox{and} \quad
P(\vecS \le {\bf c}) \le P(\vecR \le {\bf c}). \ena
\end{theorem}
\noindent {\bf Proof:} Since (\ref{PQDinq}) holds for ${\bf
U},{\bf V}$ if and only if it holds for $-{\bf U},-{\bf V}$, by
replacing the collection $\{f_\alpha({\bf x})\}_{\alpha=1}^m$ by
$\{-f_\alpha({\bf x})\}_{\alpha=1}^m$ when the functions are
decreasing, it suffices to consider the increasing case.

For $k \in \{0,\ldots,n\}$ let $ H_\alpha^k=H_\alpha \cap
\{0,\ldots,k\}, $ and with $H_\alpha$ replaced by $H_\alpha^k$, let
${\bf Z}_\alpha^k$ and $\vecS^k$ be defined as in (\ref{defYj}) and
(\ref{defSk}) respectively. We prove the first inequality in
(\ref{PQDinq}) by induction on $k$ in \bea
\label{induct-on-this-pqd} P(\vecS^k \ge {\bf c}) \le P(\vecR \ge
{\bf c}); \ena the second inequality in (\ref{PQDinq}) follows in
the same manner. Inequality (\ref{induct-on-this-pqd}) is trivially
true, with equality, when $k=0$, since then $H_\alpha^k=\emptyset$
and  ${\bf Z}_\alpha={\bf X}$ for all $\alpha \in {\bf m}$. Now
assume inequality (\ref{induct-on-this-pqd}) is true for $0 \le k <
n$ and set
$$ B=\{\alpha: k+1 \in H_\alpha\}. $$ Then \beas &&P(\vecS^{k+1} \ge {\bf
c})
\\ &=&P(f_1({\bf Z}_1^{k+1}) \ge c_1, \ldots, f_m({\bf Z}_m^{k+1})
\ge c_m) \\ &=&E[P(f_1({\bf Z}_1^{k+1}) \ge c_1, \ldots, f_m({\bf
Z}_m^{k+1}) \ge c_m|X_l,Y_l, l \not = k+1)]\\ &=&
E[P(f_\alpha({\bf Z}_\alpha^k) \ge c_\alpha, \alpha \not \in
B|X_l,Y_l, l \not = k+1) P(f_\alpha({\bf
Z}_\alpha^{k+1}) \ge c_\alpha, \alpha \in B|X_l,Y_l, l \not = k+1)]\\
&=&E[P(f_\alpha({\bf Z}_\alpha^k) \ge c_\alpha, \alpha \not \in
B|X_l,Y_l, l \not = k+1)
P(f_\alpha({\bf Z}_\alpha^k) \ge c_\alpha, \alpha \in B|X_l,Y_l, l \not = k+1)]\\
&\le&E[P(f_1({\bf Z}_1^k) \ge c_1, \ldots, f_m({\bf Z}_m^k) \ge
c_m)|X_l, Y_l,l \not = k+1]\\ &=& P(\vecS^k \ge {\bf c}) \\ &\le&
P(\vecR \ge {\bf c}), \enas where the third equality follows from
the independence of $X_{k+1}$ and $Y_{k+1}$ and the fourth from the
fact that $\{f_\alpha({\bf Z}_\alpha^k)\}_{\alpha \in B}$ has the
same conditional distribution when either $X_{k+1}$ or $Y_{k+1}$
appears as the $k+1^{st}$ coordinate of the ${\bf Z}$ vector; the
first inequality follows from the fact that conditioned on $X_l,
Y_l, l \not = k+1$, the  functions $f_\alpha({\bf Z}_\alpha^k)$ are
all increasing in $X_{k+1}$  and are therefore (conditionally)
associated, and the second inequality is the induction hypothesis
(\ref{induct-on-this-pqd}). In fact, for the first inequality above
it suffices to see that the product of the two probabilities
conditioned on $X_l, Y_l, l \not = k+1$  is the  product of
(conditional) expectations of two increasing functions of $X_{k+1}$,
which is smaller than the (conditional) expectation of the product.
$\qed$

Taking ${\bf c}=(c,\dots,c)$ we immediately  have
\begin{corollary} For all $c \in  {\bf R}$,
\label{pqd-max}
$$
P(\max_\alpha f_\alpha({\bf Z}_\alpha) \le c) \le P(\max_\alpha
f_\alpha({\bf X}) \le c) \quad {\rm or \,\,\, equivalently} \quad
\max_\alpha f_\alpha({\bf X}){\le}_{{ST}} \max_\alpha
f_\alpha({\bf Z}_\alpha).
$$
\end{corollary}
\noindent{\bf Application 1.} Consider the framework of Theorem
\ref{thm1}, with $f_\alpha({\bf x}), g_\beta({\bf x}),\alpha \in
{\cal A}, \beta \in {\cal B}$ all increasing or all decreasing
functions which depend on coordinates $K_\alpha, L_\beta$. Define
${\cal D}$ to be a collection of functions
$$
{\cal D}=\{f_\alpha({\bf X}) + g_\beta({\bf X}): K_\alpha \cap
L_\beta = \emptyset \},
$$
and for ${\bf Y}=(Y_1,\ldots,Y_n)$ as above, set
$$
{\cal D}^*=\{f_\alpha({\bf X}) + g_\beta({\bf Y}): K_\alpha \cap
L_\beta = \emptyset \}.
$$
By Theorem \ref{mon} we have
$$
{\cal D}^* \le_{PQD} {\cal D}.
$$
Applying Corollary \ref{pqd-max},
\beas
\max_{\alpha \cap \beta = \emptyset}\{f_\alpha({\bf
X}) + g_\beta({\bf X}) \} \le_{ST} \max_{\alpha \cap \beta =
\emptyset} \{f_\alpha({\bf X}) + g_\beta({\bf Y}) \}.
\enas
Exponentiating the last relation  and replacing $e^{f_\alpha}$ by
$f_\alpha$, using obvious properties of the max, we obtain
\begin{equation}
\label{lo} \max_{\alpha \cap \beta = \emptyset}\{f_\alpha({\bf X})
g_\beta({\bf X}) \} \le_{ST} \max_{\alpha \cap \beta = \emptyset}
\{f_\alpha({\bf X})g_\beta({\bf Y}) \}.
\end{equation}
and therefore
\beas
E\{\max_{\alpha \cap \beta=\emptyset} f_\alpha({\bf
X})g_\beta({\bf X})\} \le E\{\max_{\alpha \cap \beta = \emptyset}
f_\alpha({\bf X})g_\beta({\bf Y}) \} \leq E\{\max_\alpha f({\bf
X})\}\,E\{\max_\beta g({\bf X})\}, \enas for nonnegative monotone
functions $f_\alpha$ and $g_\beta$. Thus the relation (\ref{lo})
is stronger than the BKR inequality for monotone sets, which was
proved in [\ref{vbk}]. Alexander [\ref{alex}] presents similar
functional versions in this context.

As an example we return to order statistics as in Section
\ref{order-statistics}. From (\ref{lo}) we derive, for example, that
$$X_{[1]}X_{[2]} \le_{ST} X_{[1]}Y_{[2]} \vee Y_{[1]}X_{[2]}.$$
Generalizing by using the functions (\ref{fgprod}), we obtain for
any $p+q=m$,
$$
\prod_{j=1}^m X_{[j]} \le_{ST} \max_{\{i_1,\ldots,i_p\} \cup
\{j_1,\ldots,j_q\} =\{1,\ldots,m\}}\prod X_{[i_q]}Y_{[j_q]}.
$$

\section{Appendix on Measurability}
\label{appendix} In this section we briefly deal with various
measurability issues. The measurability of the functions defined in
(\ref{Kalphax}) can be seen from
$$
{\tilde f}_{K({\bf x})}({\bf x})=\sum_K {\tilde f}_K({\bf x}){\bf
1}(K({\bf x})=K),
$$
since the given function $K({\bf x})$ is assumed measurable.
Similarly for (\ref{abuse}), $$ {\tilde f}_{\Gamma({\bf x})}({\bf
x})=\sum_{A \in 2^{\cal A}} \sup_{\alpha \in A}f_\alpha({\bf
x}){\bf 1}(\Gamma({\bf x})=A).$$

We next prove that given a non-negative, $(\bbs,\mathbb{B})$
measurable function $f:\spa \rightarrow {\bf R}$ and any $K
\subseteq \N$, the function $\underline{f}_K({\bf x})$ defined in
(\ref{fgismin}) is $(\bbs,\mathbb{B})$ measurable. Letting
$$
f_r({\bf x})=\min(f({\bf x}),r)
$$
and $\PP_L({\bf x})$ be the marginal of $\PP({\bf x})$ on the
coordinates ${\bf x}_L$, we have \beas \lim_{p \rightarrow \infty}
\left( \int (r-f_r({\bf x}))^p d\PP_{K^c}({\bf x})\right)^{1/p} =
{\rm ess} \sup_{{\bf y} \in [{\bf x}]_K} (r-f_r({\bf y}))= r -
{\rm ess} \inf_{{\bf y} \in [{\bf x}]_K} f_r({\bf y}). \enas
Tonelli's theorem (see e.g.~[\ref{folland}]) now implies that
${\rm ess} \inf_{{\bf y} \in [{\bf x}]_K} f_r({\bf y})$ is
measurable. Letting $r \uparrow \infty$ shows that (\ref{fgismin})
is measurable.

The only complication regarding measurability of the pair $(K({\bf
x}),L({\bf x}))$ in (\ref{maxKL-measurable}) is that the maximum
may not be uniquely attained, since otherwise we would simply have
$$ \{{\bf x}: K({\bf x})=K, L({\bf x})=L
\} = \bigcap_{K' \cap L' = \emptyset }\{ {\bf x}: {\underline
f}_K({\bf x}){\underline g}_L({\bf x}) \ge {\underline
f}_{K'}({\bf x}){\underline g}_{L'}({\bf x})\},
$$ a finite intersection of measurable sets, so measurable. To
handle the problem of non-uniqueness, let $\prec$ be an arbitrary
total order on the finite collection of subsets of $\N \times \N$,
so that when the max is not unique we can choose $(K({\bf
x}),L({\bf x}))$ to be the first disjoint pair that attains the
maximum. Then $ \{{\bf x}: K({\bf x})=K, L({\bf x})=L \} = F \cap
G$ where
$$F=\bigcap_{\stackrel{(K',L') \prec  (K,L)}{K' \cap L' = \emptyset} }
\{ {\bf x}: {\underline f}_K({\bf x}){\underline g}_L({\bf x}) >
{\underline f}_{K'}({\bf x}){\underline g}_{L'}({\bf x})\}$$ and
$$G=\bigcap_{\stackrel{(K',L')  \succeq (K,L)}{K' \cap L' = \emptyset} }\{ {\bf x}: {\underline
f}_K({\bf x}){\underline g}_L({\bf x}) \ge {\underline
f}_{K'}({\bf x}){\underline g}_{L'}({\bf x})\}$$ and again
measurability follows. \lcolor{Similar remarks apply to the maximizing $K({\bf x}),L({\bf x})$ in \eqref{I-to-II} for the implication $\dBKR \Rightarrow  \aBKR$.}

\lcolor{Finally, we note inequality \aBKR, and therefore also \bBKR, holds on the completion of $(\spa,\overline{\bbs})$ of $(\spa,\bbs)$ with respect to $P$. Proposition 2.12 of [\ref{folland}] shows that for every $\overline{\bbs}$ measuable function $\overline{f}$ there exists an $\bbs$ measurable function $f$ such that $\overline{f}=f$ with (completed) measure one. Hence, replacing all $\overline{\bbs}$ measurable functions in \aBKR\, by their $\bbs$ measurable counterparts and applying \aBKR\, over the space $(S,\bbs)$ shows \aBKR\, holds on $(S,\overline{\bbs)}$.}\\[1ex]

\noindent Acknowledgment:
The authors would like to thank Richard Arratia for pointing out to us that \eqref{was.5} is an inequality in general, and for other insightful comments.

\section*{Bibliography}

\begin{enumerate}
\item \label{alex} Alexander, K. (1993) A note on some rates of
convergence in first-passage percolation. {\em Ann. Appl. Probab.}
{\bf 3}, 81--90.

\item \label{AGH} Arratia, R., Garibaldi, S. and Hales, A. (2015) The van den Berg-Kesten-Reimer inequality for
infinite spaces. Unpublished manuscript.

\item \label{epw} Esary, J. D., Proschan, F., Walkup, D. W. (1967)
Association of random variables, with applications. {\em Ann.
Math. Statist.} {\bf 38},  1466--1474.

\item \label{folland} Folland, G. (1999) Real Analysis: Modern
Techniques and Their Applications, Wiley, N.Y.

\item \label{GoRi} Goldstein, L. and Rinott, Y. (2007) Functional BKR Inequalities, and their Duals, with Applications. {\em J. Theor. Probab.}, {\bf 20}, 275-293.

\item \label{kss} Kahn, J., Saks, M., and Smyth, C. (2000) A dual
version of Reimer's inequality and a proof of Rudich's conjecture.
15th Annual IEEE Conference on Computational Complexity 98-103,
IEEE Computer Soc., Los Alamitos, CA, 2000.

\item \label{kr} Karlin, S. and Rinott Y (1980). Classes of
orderings of measures and related correlation inequalities. I.
Multivariate totally positive distributions. {\em J. Multivariate
Anal.} {\bf 10}, 467--498.

\item \label{reim} Reimer, D. (2000) Proof of the van den
Berg-Kesten conjecture.  {\em Combin. Probab. Comput.}  {\bf 9},
27--32.

\item \label{ss} Shaked, M., and Shanthikumar, J. (1994).
Stochastic orders and their applications. Probability and
Mathematical Statistics. Academic Press, Inc., Boston, MA.

\item \label{vbf} van den Berg, J., and Fiebig, U. (1987) On a
Combinatorial Conjecture Concerning Disjoint Occurrences of Events
{\em Ann. Prob.} {\bf 15}, 354-374.

\item \label{vbk} van den Berg, J., and Kesten, H. (1985)
Inequalities with applications to percolation and reliability.
{\em J. Appl. Probab.}  {\bf 22}, 556--569.

\end{enumerate}

\end{document}